\documentclass[12pt, a4paper]{article}

\usepackage[latin1]{inputenc}
\usepackage[english]{babel}
\usepackage{mathtext,amssymb,amsmath}
\usepackage{amsthm,amsfonts,mathrsfs}
\usepackage{color,graphicx}
\usepackage[all]{xy}
\usepackage{psfrag}
\usepackage{float}

\author{In\^es Silva de Oliveira and Paul A. Schweitzer, S.J.}
\title{Locally convex hypersurfaces immersed in $H^n \times R$}

\newtheorem{theo}{Theorem}[section]

\newtheorem{prop}[theo]{Proposition}

\newtheorem{lemma}[theo]{Lemma}
\newtheorem*{maintheo}{Main Theorem}
\newtheorem{hyp}[theo]{Hypothesis}

\newcommand{\Hset}{\mathbb{H}}
\newcommand{\Rset}{\mathbb{R}}
\newcommand{\Bset}{\mathbb{B}}
\newcommand{\Sset}{\mathbb{S}}

\newcommand{\HxR}{\Hset^n\times\Rset}

\newcommand{\HkxR}{\Hset^k\times\Rset}
\newcommand{\Hkxo}{\Hset^k\times\{0\}}
\newcommand{\Sinf}{\Sset^{n-1}_\infty}
\newcommand{\Sinfk}{\Sset^{k-1}_\infty}

\begin{document}

\maketitle

\begin{abstract}
We prove a theorem of Hadamard-Stoker type: a connected locally convex complete hypersurface immersed in $\HxR$ ($n\geq 2$), where $\Hset^n$ is $n$-dimensional hyperbolic space, is embedded and homeomorphic either to the $n$-sphere or to $\Rset^n$. In the latter case it is either a vertical graph over a convex domain in $\Hset^n$ or has what we call a simple end.

\medskip

\noindent{\bf keywords:} Hadamard-Stoker theorems, convex hypersurface, positive extrinsic curvature, hyperbolic geometry

\medskip


\end{abstract}



\section{Introduction}

Let $\Sigma$ be a complete hypersurface immersed in $\HxR$, the product of $n$-dimensional hyperbolic space $\Hset^n$ with the line, and let $\pi:\HxR\to \Hset^n$ be the projection on the first factor. A point $\theta$ in the sphere at infinity $\Sinf$ of $\Hset^n$ is a {\bf simple end} of $\Sigma$ if $\theta$ is the unique point of accumulation of $\pi(\Sigma)$ in $\Sinf$ and
every complete totally geodesic hyperplane $Q\subset \Hset^n$ either has $\theta$ as an accumulation point or else meets $\pi(\Sigma)$ in a compact set. We prove the following.

\begin{maintheo} Let $\Sigma$ be a complete, connected hypersurface immersed in $\Hset^n\times\Rset$ with positive definite second fundamental form, $n\geq 2$. Then $\Sigma$ is embedded and homeomorphic either to the $n$-sphere or to Euclidean space $\Rset^n$. In the latter case, it is either a vertical graph over an open convex set or has a simple end. When $n=1$ the theorem holds under the additional hypothesis that $\Sigma$ is embedded.
\label{mainthm} \end{maintheo}

The case $n=2$ was proven by J. Espinar, J. G\'alvez, and H. Rosenberg \cite{EGR}, so our paper is an extension of their work to higher dimensions. We prove the theorem by induction on $n$, beginning with the easy case $n=1$, when $\Sigma$ is assumed to be embedded. Our paper reproves the case $n=2$; we use many ideas of \cite{EGR}, but our proof does not depend on their result.

\bigskip
\noindent{\bf Historical Background.}
In 1897, J. Hadamard \cite{H} proved a theorem about compact, locally strictly
convex surfaces in the Euclidean space $\mathbb{R}^3$, showing that
such surfaces are embedded and homeomorphic to the sphere. Since
then many generalizations have adapted the assumptions about
the curvature and considered new spaces in which these surfaces
could be immersed in order to obtain analogous results.

Important contributions were made by J. Stoker \cite{S},
S.S. Chern and R. K. Lashof \cite{CL}, R. Sacksteder \cite{Sa}, M. do Carmo and E.
Lima \cite{ME}, M. do Carmo and F. Warner \cite{MW}, R.J.Currier \cite{C}, S. Alexander
\cite{A}, and I. Tribuzy \cite{I}.

In 2009, J. Espinar, J. Gálvez and H. Rosenberg \cite{EGR} extended the Hada\-mard-Stoker theorem for surfaces immersed in $\mathbb{H}^2
\times \mathbb{R}$ assuming that such a surface is connected and
complete with all principal curvatures positive. They proved that such a surface is properly embedded and homeomorphic to the
sphere if it is closed or to the plane $\mathbb{R}^2 $ if not. In the second case, $\Sigma$ is a graph over a convex domain in $\mathbb{H}^2 \times
\{0\} $ or $\Sigma$ has a simple end.

J. Espinar and H. Rosenberg \cite{ER} proved in 2010 that if $ \Sigma$ is a locally strictly
convex, connected hypersurface properly immersed in $M ^ n \times \mathbb{R} $, where $ M^n $ is a $
\emph{1/4-pinched}$ manifold, then $ \Sigma $ is properly embedded and
homeomorphic to $\mathbb{S}^n$ or $ \mathbb{R}^n$. In the second case, $
\Sigma $ has a $ \emph{top} $ or $\emph{bottom} $ end, where $\Sigma $
has a $\emph {top} $ (respectively, \emph {bottom}) end $E$ if for any divergent sequence $ \{p_{n}\} \subset E $ the height function  $h: \Sigma \rightarrow E $ goes to $+ \infty$ (respectively $-
\infty$).

Also in 2010, J. Espinar and the first author \cite{EO} showed a
Hadamard-Stoker-Currier type theorem for surfaces immersed in a
$3$-dimen\-sional Riemannian
manifold $\mathcal{M}(\kappa, \tau)$ that fibers over a $2$-dimensional Riemannian
manifold $\mathbb{M}^2$ so that the fibers are the trajectories of a unit Killing vector
field. More precisely, $\mathbb{M}^2$ is required to be a strict
Hadamard surface, i.e., $ \mathbb{M}^2 $ has Gaussian curvature $\kappa$ less than a negative constant, and $\tau$ is the curvature of the bundle.

\bigskip
\noindent{\bf Open Questions.}
There are other spaces to be studied, for example: $ \mathbb{S}^n
\times \mathbb{R}$, $ \mathbb{M} \times \mathbb{R} $, where
$\mathbb{M} $ is a Hadamard manifold of dimension $n$, Heisenberg
spaces of dimension $n$, among others. It possible that results of Hadamard-Stoker type are valid in these spaces.

\bigskip

Preliminary results and the case $n=1$ are treated in Section \ref{preliminary}. The proof of the Main Theorem in the case of a vertical graph is given in Section \ref{vgraph} and the remaining cases of the theorem are proven in Section \ref{othercases}. Examples of hypersurfaces with simple ends are given in Section \ref{example}.

This paper is based on the doctoral thesis \cite{IO} of the first author at
the Pontifical Catholic University of Rio de Janeiro under the
direction of the second author.


\section{Preliminary results and the case $n=1$}\label{preliminary}

Let $\Hset^k$ be  the $k$-dimensional hyperbolic space, with a Riemannian metric of constant curvature $-1$ (See \cite{M}, \cite{ST} for details of Riemannian geometry). In the Poincar\'e disk model on the open unit ball $\Bset^k=\{(x_1,\dots,x_k)\in \Rset^k\ |\ |x|<1\}$, where $|x|$ is the Euclidean norm given by $|x|^2=\sum_{i=1}^k x_i^2$,
the Riemannian metric is defined to be
$$ds^2= \frac{4\sum_{i=1}^k dx_i^2}{(1-|x|^2)^2}.$$
In this model the sphere at infinity, $\Sinfk$, is
the unit sphere in $\Rset^k$, so that $\Hset^k\cup \Sinfk$ is a compact $k$-dimensional ball.
If $Q$ is a totally geodesic hyperplane in $\Hset^k$, we say that $P=Q\times \Rset$, which is also totally geodesic, is a {\bf vertical hyperplane} in $\HkxR$. If $\gamma$ is a complete geodesic in $\Hkxo$ parametrized by $t\in \Rset$, we let
$P^{k}_{\gamma}(t)$ be the vertical hyperplane orthogonal to $\gamma$ at the point $\gamma(t)$, and then these vertical hyperplanes $P^{k}_{\gamma}(t)$ form a foliation of $\HkxR$. Let $\pi: \HkxR\to \Hset^k$ denote the projection on the first factor.


The following proposition will be used in the inductive step of the proof of the Main Theorem when $M=\HkxR$, but we state it in a more general context.

\begin{prop} Let $M$ be a $(k+1)-$dimensional Riemannian manifold and $\Sigma^k$ a hypersurface immersed in $M$ with strictly positive second fundamental form, $k\geq 2$, and let $P$ be a totally geodesic hypersurface in $M$. If $\Sigma^k$ and $P$ intersect
transversally, then every connected component $\Sigma^{k-1}$ of
$\Sigma^k \cap P$ is a $(k-1)-$dimensional hypersurface in $P$ with second fundamental form $II(\Sigma^{k-1})>0$.

\end{prop}

\noindent{\bf Proof.}
Let $\gamma(t)$ be a curve in $\Sigma^{k-1}$, where $t$ is
parametrized by arc length. As $P$ is totally geodesic,
we have $\overline{\nabla}_{\gamma'}(\gamma')=
\nabla_{\gamma'}^{P}(\gamma')$, where $\nabla^{P}$ and
$\overline{\nabla}$ are the connections in $P$ and
$M$ respectively. Let $N_{\Sigma^{k-1},P}$ and $N_{\Sigma^{k}}$ be the unit normal vectors to $\Sigma^{k-1}$ in $P$ and to $\Sigma^{k}$ in $M$, respectively. We want
$\langle \overline{\nabla}_{\gamma'}(\gamma'),
N_{\Sigma^{k-1},P}\rangle$ to be positive. Note that  $\langle
\overline{\nabla}_{\gamma'}(\gamma'), N_{\Sigma^{k}} \rangle$ is
positive, because by hypothesis  the second fundamental form of
$\Sigma^n$ is strictly positive.

Writing $N_{\Sigma^{k}}(\gamma(t)) = \overline{N}+N^{\bot}$, where
$\overline{N}$ is tangent to $P$ and $N^{\bot}$ is orthogonal to $P$, and taking the inner product with
$\overline{\nabla}_{\gamma'}(\gamma')$, we obtain

\begin{align}
0<\langle \overline{\nabla}_{\gamma'}(\gamma'), N_{\Sigma^{k}} \rangle
& = \langle \overline{\nabla}_{\gamma'}(\gamma'), \overline{N} +
N^{\bot}\rangle \nonumber \\ & = \langle
\overline{\nabla}_{\gamma'}(\gamma') ,\overline{N}\rangle. \nonumber
\end{align}
 Note that $\overline{N}$ is orthogonal to $\Sigma^{k-1}$ and tangent to $P$. Hence $\overline{N}$ is a multiple of
 ${N}_{\Sigma^{k-1},P}$, so we can write
\begin{align}
 \langle \nabla_{\gamma'}^{P}(\gamma'),\overline{N}\rangle & = \langle
 \nabla_{\gamma'}^{P}(\gamma'),|\overline{N}| N_{\Sigma^{k-1},P}\rangle \nonumber \\ & = |\overline{N}|\langle \nabla_{\gamma'}^{P}(\gamma'),
 N_{\Sigma^{k-1},P}\rangle, \nonumber
\end{align}
and therefore $\langle \nabla_{\gamma'}^{P}(\gamma'),
 N_{\Sigma^{k-1},P}\rangle > 0$, as claimed.  \qed

 \vspace{0.5cm}


Next we begin the proof of the Main Theorem. The proof is by induction on the dimension $n$ of $\Sigma^n$. The proof for the case $n=1$ is given in this section, and the proof for $n=k>1$, when the
theorem is supposed true for $n=k-1$, is given in the next two sections.

\bigskip
\noindent{\bf The Proof for $n=1$.} When $n=1$, we assume the additional hypothesis that the curve $\Sigma^1$ is embedded in $\Hset^1\times \Rset$, which is isometric to $\Rset^2$. If $\Sigma^1$ is not a vertical graph over an open interval in $\Hset^1$, which is one of the possible conclusions of the theorem, then there must be a point
$p_0\in\Sigma^1$ where the tangent line is vertical. If there is a second such point, then $\Sigma^1$ closes up to a circle, and if $p_0$ is the only such point, then both ends of $\pi(\Sigma^1)$
must converge to one of the two points at infinity of $\Hset^1$, which gives a simple end.
This completes the proof of the theorem for $n=1$. \qed


\section{The case of a vertical graph}\label{vgraph}

In the previous section the Main Theorem was proven for the case $n=1$. We complete the proof in this section and the next one. In both sections we suppose that the theorem holds for $n=k-1\geq 1$ and we assume the hypotheses of the Main Theorem for the case $n=k$ in order to prove the inductive step.
Thus we assume that $\Sigma^k$ is a
complete connected $k$-dimensional manifold immersed in
$\HkxR$, which has the product Riemannian metric, and we suppose that the second fundamental form of $\Sigma^k$ is positive definite.

For the rest of this section we also assume the following.

\begin{hyp} No vertical hyperplane is tangent to $\Sigma^k$.\label{hyp1}
\end{hyp}

We shall prove the Main Theorem under this additional Hypothesis by showing that in this case $\Sigma^k$ is the graph of a smooth function $f: \Omega\to \Rset$ where $\Omega$ is an open convex domain in $\Hset^k$. Note that the Hypothesis implies that every vertical hyperplane $P$ is transverse to $\Sigma^k$.

\begin{lemma} Hypothesis \ref{hyp1} implies that if $P$ is a vertical hyperplane in $\HkxR$, then each connected component $\Sigma^{k-1}$ of the intersection $\Sigma^k\cap P$ is a graph over a convex domain in $\pi(P)$. \label{verticalgraph}
\end{lemma}

\noindent{\bf Proof.} By Proposition 2.1, $\Sigma^{k-1}$ is a $(k-1)-$dimensional surface with strictly positive second fundamental form in
$P$. Hence we can apply the induction hypothesis, so $\Sigma^{k-1}$ must be
homeomorphic to $\mathbb{S}^{k-1}$ or $\mathbb{R}^{k-1}$, and in the latter case it must have a simple end or be a vertical graph over a convex domain.

Now if $\Sigma^{k-1}$ is homeomorphic to $\mathbb {S}^{k-1}$, then there exists a point in $\mathbb{S}^{k-1}$ with a
$(k-1)$-dimensional vertical tangent plane, which means that $\Sigma^k$ has
a vertical tangent $k$-plane at this point, contradicting Hypothesis \ref{hyp1}, so this case is excluded.
If $\Sigma^{k-1}$ is homeomorphic to $\mathbb{R}^{k-1}$ and has a simple end $\theta$ in $\Sset^{k-1}_\infty$, the sphere at infinity of $\pi(P)$, let $\beta(t)$ be a complete horizontal geodesic in $\pi(P)\times\{0\}\subset P$ which converges to the point
$(\theta,0)$ as $t\to\infty$. Consider the foliation of $P$ by vertical $(k-1)$-hyperplanes $P^{k-1}_\beta(t)$ orthogonal to $\beta$, where $P^{k-1}_\beta(t)$ meets $\beta$ at $\beta(t)$. If $\bar t$ is the smallest value of $t$ such that $\Sigma^{k-1}\cap P^{k-1}_\beta(t)$ is non-empty, then the vertical hyperplane $P^{k-1}_\beta(\bar t)$
will be tangent to $\Sigma^{k-1}$, so $\Sigma^k$ will also have a
tangent vertical hyperplane, contradicting Hypothesis \ref{hyp1}.
Thus $\Sigma^{k-1}$ must be a vertical graph over a convex domain in $\pi(P)$.
\qed

\begin{lemma}
Let $\{P^{k}_{\gamma}(t)\}$ be the foliation of $\HkxR$ by vertical hyperplanes orthogonal to a complete geodesic $\gamma$ in $\Hkxo$ and let $\Sigma^{k-1}(0)$ be a component of
$\Sigma^{k} \cap P_{\gamma}(0)$.
Let $\Sigma^{k-1}(t)\subset \Sigma^{k} \cap P_{\gamma}(t)$ be the continuous variation of $\Sigma^{k-1}(0)$
as $t$ varies. Then $\Sigma^{k-1}(t)$ cannot become disconnected for any $t$.\label{connected}
\end{lemma}

\noindent{\bf Proof.} Note that $\Sigma^{k-1}(t)$ is well defined by transversality.
By Lemma \ref{verticalgraph}, every component of $\Sigma^{k-1}(t)$ is the graph of a function defined on a convex open domain in $\pi(P_{\gamma}(0))$, which is isometric to
$\Hset^{k-1}$. By transversality the set of $t$ with
$\Sigma^{k-1}(t)$ connected is open.
Suppose, to find a contradiction, that $\Sigma^{k-1}(t)$ is not connected for some $t$, say $t>0$. Then there will be a smallest positive $\bar t$ such that $\Sigma^{k-1}(\bar t)$ is not connected. Take $x_{\bar t}$ and $y_{\bar t}$ to be points in distinct components of $\Sigma^{k-1}(\bar t)$. By transversality there exist $\delta>0$ and continuous curves $t\mapsto x_t$ and $t\mapsto y_t$ defined for $t\in [\bar t-\delta,\bar t]$ with $x_t, y_t\in\Sigma^{k-1}(t)$. For
$t\in [\bar t-\delta,\bar t)$, $\Sigma^{k-1}(t)$ is connected
and thus a graph over a convex open domain in
$\pi(P_{\gamma}(t))$. Let $\alpha_t$ be the geodesic in
$\pi(P_{\gamma}(t))$ joining $\pi(x_t)$ and $\pi(y_t)$ and take $A_t=\pi^{-1}(\alpha_t)\cap \Sigma^{k-1}(t)$, which is a graph over $\alpha_t$ and a curve joining $x_t$ to $y_t$ in $\Sigma^{k-1}(t)$. The limit of $A_t$ as $t$ tends to $\bar t$
is a curve $A_{\bar t}$ in $\Sigma^{k-1}(\bar t)$, which is complete, since $\Sigma^k$ is complete, but $A_{\bar t}$ cannot be connected since it contains
$x_{\bar t}$ and $y_{\bar t}$ which are in different components. It follows that $A_{\bar t}$
must diverge vertically to $-\infty$ or $+\infty$, but since $A_{\bar t}$ has strictly positive curvature that is impossible. \qed

\bigskip

In view of the previous lemma, the union $\cup\Sigma^{k-1}(t)$ is an open and closed set in $\Sigma^k$, so it must be the whole connected set $\Sigma^k$. Consequently $\Sigma^k$ is a vertical graph over a set $\Omega$ in $\Hset^k$. Since $\Omega$ is diffeomorphic to the $k$-dimensional manifold $\Sigma^k$ under the projection $\pi$, it is open in $\Hset^k$. Now given any two points $x,y\in\Omega$, let $P$ be a vertical hyperplane containing both $(x,0)$ and $(y,0)$, and apply the argument of Lemma \ref{connected} to see that $\Sigma^k\cap P$ must be connected.
Then $\Sigma^k\cap P$ is a vertical graph over a convex domain $\Omega_P$ in $\pi(P)$. The geodesic from $x$ to $y$ lies in $\Omega_P\subset \Omega$, so $\Omega$ is convex. Now the convex open set $\Omega$ in $\Hset^k$ is homeomorphic to $\Rset^k$, so $\Sigma^k$ is also homeomorphic to $\Rset^k$.  This completes the proof of the Main Theorem under Hypothesis \ref{hyp1}.


\section{The proof for the remaining cases}\label{othercases}

In this section we complete the proof of the Main Theorem by proving it in the case that the previous Hypothesis \ref{hyp1} does not hold, i.e., under the following assumption.

\begin{hyp} There is a vertical hyperplane $P_0$ that is tangent to $\Sigma^k$ at a point $p_0\in\Sigma^k$.\label{hyp2}
\end{hyp}

\begin{figure}[htbp]
\centering
\includegraphics[width=5.0cm]{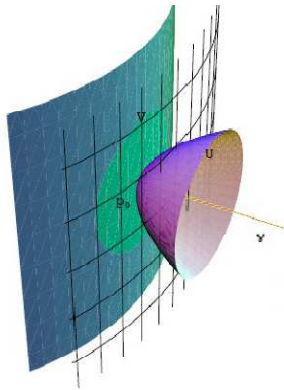}
\caption{Near $p_0$, $\Sigma^k$ lies on the external side of $P_\gamma(\overline{t})$.} \label{figviz}
\end{figure}

Without loss of generality, we suppose that $p_0\in\Hkxo$.
In view of the hypothesis that the hypersurface $\Sigma^k$ has strictly positive curvature, there is a neighborhood $\mathcal{U}_0$ of $p_0$ in
$\Sigma^k$ that lies entirely on one side of $P_0$, except at the point $p_0$; we shall call this the external side (See Figure \ref{figviz}).
Let $\gamma$ be the geodesic in $\Hkxo$ orthogonal to $P_0$ (and therefore also to $\Sigma^k$) at $p_0$, parametrized by arclength, with $\gamma(0)=p_0$, oriented so that $\gamma(0,\infty)$ is on the external side of $P_0$. Let $P_\gamma(t)$ be the vertical $k-$plane orthogonal to $\gamma$ at $\gamma(t)$, and note that these
vertical planes for all $t\in {\mathbb R}$ form a smooth foliation of $\HkxR$.
Denote by $\Sigma^{k-1}(t)$ the connected component of $\Sigma^k\cap P_\gamma(t)$ that is the continuation of $\Sigma^{k-1}(0)= \{p_0\}\subset \Sigma^k\cap P_\gamma(0)$. We do not exclude the possibility that there may be other
components of $\Sigma^k \cap P_\gamma(t)$, but we only consider
$\Sigma^{k-1}(t)$. Note that for $t>0$ close to $0$,
$\Sigma^{k-1}(t)$ is homeomorphic to the $(k-1)-$sphere
${\mathbb S}^{k-1}$.
Furthermore, if $\Sigma^{k-1}(t)$ is compact and
$\Sigma^k$ is transverse to $P_\gamma(t)$ along $\Sigma^{k-1}(t)$
for all $t$ in an interval $(0,t_0)$, the sets
$\Sigma^{k-1}(t)$ vary continuously and are all diffeomorphic to ${\mathbb S}^{k-1}$ (See \cite{Mi}, Theorem 3.1).

Now there are four cases to be considered.

\bigskip

\noindent{\bf Case 1.} The set $\Sigma^{k-1}(t)$ is compact and
$\Sigma^k$ is transverse to $P_\gamma(t)$ along $\Sigma^{k-1}(t)$
for $0<t<\overline{t}$, but $\Sigma^k$ is not transverse to $P_\gamma(\overline{t})$ at some point $p_1\in\Sigma^{k-1}(\overline{t})$.

\vspace{0.5cm}

\begin{figure}[htbp]
\centering
\includegraphics[width=5.0cm]{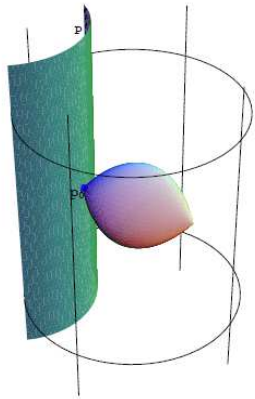}
\caption{$\Sigma^k$ homeomorphic to the $k-$sphere  $\mathbb{S}^k$.} \label{casob1b}
\end{figure}

As in the case of $p_{0}$, the point $p_1$ has
a neighborhood $\mathcal{U}_1$ in $\Sigma^k$
such that $\mathcal{U}_1 \setminus \{p_{1}\}$ lies entirely on one side of $P_{\gamma}(\overline{t})$. Since $\Sigma^{k-1}(\overline{t})$ is the continuation of $\Sigma^{k-1}(t)$ with $t<\overline{t}$, $\mathcal{U}_1\setminus\{p_1\}$ must lie on the side of $P_\gamma(\overline{t})$ with $t<\overline{t}$. The sets
$\Sigma^{k-1}(t)$ for $0<t< \overline{t}$
are all diffeomorphic to ${\mathbb S}^{k-1}$ and their union with the two points $p_0$ and $p_1$ is homeomorphic to the sphere ${\mathbb S}^{k}$. This union is open and closed in
$\Sigma^k$, so by connectedness it must coincide with
$\Sigma^k$, which is therefore a topological $k-$sphere (See Figure \ref{casob1b}). This completes the proof in Case 1.

\bigskip
In the remaining cases we exclude Case 1, so if $\Sigma^{k-1}(t)$ is compact for all $t$ in some interval $(0,t_0)$, then
 $\Sigma^k$ is transverse to the planes
$P_\gamma(t)$ at all points in these sets $\Sigma^{k-1}(t)$.

\bigskip 
\noindent{\bf Case 2.} The intersection $\Sigma^{k-1}(t)$ is compact and non-empty for all $t>0$.

\vspace{0.5cm}

\begin{figure}[htbp] 
\centering
\includegraphics[width=4.0cm]{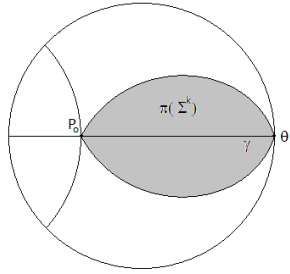}
\caption{Projection of $\Sigma^k$ on $\mathbb{H}^k$ in Case 2.}
\end{figure}

\vspace{0.5cm}

If $\Sigma^{k-1}(t) \subset\Sigma^{k}\cap P_{\gamma}(t)$ remains
compact and non-empty for all  $t>0$, then the sets
$\Sigma^{k-1}(t)$ for $t>0$ are $(k-1)-$spheres that are embedded in $P_{\gamma}(t)$ by the inductive hypothesis. Hence the union of these sets with the point $p_0$ must be all of the connected hypersurface $\Sigma^k$, which will also be embedded and
homeomorphic to $\mathbb{R}^k$. If we let
$\theta\in\mathbb{S}^{k-1}_{\infty}$ be the limit point of $\gamma(t)$ as $t\to\infty$, then
$\partial_{\infty}\pi(\Sigma^{k})=\{\theta\}$, since $\theta$ is the only point of $\mathbb{S}^{k-1}_{\infty}$ that is on the external side of all the planes $P_{\gamma}(t)$.

Furthermore, if $Q$ is any complete
totally geodesic hypersurface in $\Hset^k$ which does not have
$\theta$ as a point of accumulation, then
$\pi(\Sigma^k)\cap Q$ cannot have any point of accumulation at
infinity. Hence it must be a closed and bounded set in $Q$, so it is compact. Therefore in this case $\theta$ is a simple end
of $\Sigma^k$.


\bigskip
\noindent{\bf Case 3.} The intersection $\Sigma^{k-1}(t)$ remains compact and non-empty for $0\leq t<\bar t$ but becomes empty for every $t>\bar t$, for some $\bar t>0$.

\bigskip

Note first of all that $\Sigma^k$ cannot intersect $P_{\gamma}(\overline{t})$ transversely at any point of $\Sigma^{k-1}(\bar t)$, for otherwise $\Sigma^{k-1}(t)$ would be non-empty for values of $t$ slightly greater than $\bar t$, contrary to the hypothesis of Case 3. Hence any point in $\Sigma^{k-1}(\bar t)$ would have a vertical tangent plane, a situation which has already been treated in Case 1, so we may suppose that $\Sigma^{k-1}(\bar t)$ is empty.
Consequently $\cup_{t<\bar t} \Sigma^{k-1}(t)$ is both open and closed in
$\Sigma^k$, so it must be all of $\Sigma^k$. As in Case 2, we see that $\Sigma^k$, which is the union of the point $p_0$ and the topological $(k-1)-$spheres $\Sigma^{k-1}(t)$ for $0<t<\bar t$, is homeomorphic to ${\mathbb R}^k$.

Now take a sequence of points $p_n\in \Sigma^{k-1}(t_n)$ such that
$t_n\to \bar t$. Since $\overline \Hset^k =\Hset^k\cup\partial\Hset^k$
is compact, there must be a subsequence of $\{\pi(p_n)\}$ that
converges to a point $\theta$ in $\overline \Hset^k$. Since we are
supposing that $\Sigma^{k-1}(\bar t)$ is empty, $\theta$ must be
in ${\mathbb S}^{k-2}_\infty(\bar t) = \partial\pi(P_\gamma(\bar t))$
and must be
an accumulation point of $\pi(\Sigma^k)$ at infinity.
To complete the proof in this case, we must show that $\theta$ is the only accumulation point and that it is a simple end.

Suppose there were another accumulation point of $\pi(\Sigma^k)$, say $\widetilde{\theta}\in \Sinfk$.
Then $\widetilde{\theta}$, like $\theta$ and any other accumulation
points of $\pi(\Sigma^k)$, must be in ${\mathbb S}^{k-2}_\infty (\bar t)$,
since there are no accumulation points for $t<\bar t$ and $\Sigma^k\cap
P_\gamma(t)$ is empty for $t>\bar t$. Choose the parametrization
of $\mathbb{S}^{k-2}_{\infty}(\bar t)$ in the disk model of $\Hset^k$ so that $\widetilde{\theta}$ is the antipode of $\theta$.
Let $\mu$ be the totally geodesic $2-$plane in $\Hkxo$ that contains both $p_{0}$ and the geodesic from $(\theta,0)$ to $(\widetilde{\theta},0)$. Take
$\xi\subset \HkxR$ to be the totally geodesic vertical
$(k-1)-$plane orthogonal to $\mu$
that meets $\mu$ at the single point $p_0$.
Consider the $1-$parameter family of vertical $k-$planes $\{P(\varphi)=Q(\varphi)\times{\mathbb R}\}$ containing
$\xi$ for $ \varphi \in [-\pi/2,3\pi/2] $, parametrized injectively so
that $\theta$ is a point at infinity of $Q(0)$, $\widetilde {\theta}$ is a point at infinity of $Q(\pi)$, and the
tangent plane at $p_0$ is not among the planes $P(\varphi)$. Now the intersection
$$\Sigma^{k-1}(\varphi)= \Sigma^k\cap P(\varphi)$$ is a $(k-1)-$surface in $P(\varphi)$ which satisfies
the hypotheses of the Main Theorem, so by the inductive hypothesis it must be homeomorphic to
$\mathbb{S}^{k-1}$ or $ \mathbb{R}^{k-1}$ and in the latter case $
\Sigma^{k-1}(\varphi)$ it is either a vertical graph over a convex domain
in $ \mathbb{H}^{k-1}$ or it has a simple end. However,
$\Sigma^{k-1}(\varphi) $ can not be a vertical graph, because this
$(k-1)-$surface contains the point $p_{0}$ where there
is a vertical tangent hyperplane.
Moreover, for each parameter $\varphi$
with $\varphi<$ 0 or $\varphi >
\pi$ the intersection $\Sigma^{k-1}(\varphi)$ is a bounded complete $(k-1)-$surface contained in
$\cup_{0 \leq t \leq t_0}
\Sigma^{k-1}(t)$ for some $t_0<\overline{t}$ and is therefore compact, with no accumulation points at infinity.

If for some $\overline{\varphi}$ with $0 < \overline{\varphi} < \pi,
\Sigma^{k-1}(\overline {\varphi}) $ is compact, then by the inductive hypothesis it will be
homeomorphic to the sphere $\mathbb{S}^{k-1}$ and must bound
a closed $k-$dimensional ball $D$ in $\Sigma^k$, which is homeomorphic to $\Rset^k$. The ball $D$ must coincide either with $\Sigma^k\cap \cup_{\varphi\leq\overline {\varphi}} Q(\overline {\varphi})$
or else with $\Sigma^k\cap \cup_{\varphi\geq\overline {\varphi}} Q(\overline {\varphi})$; in the first case, $\theta$ cannot be a limit point of $\Sigma^k$, and in the second case, $\widetilde{\theta}$ cannot be a limit point of $\Sigma^k$, contradicting the choice of $\theta$ and $\widetilde{\theta}$. Hence for $0 <\varphi < \pi, \Sigma^{k-1}(\varphi) $  must be noncompact and it must have a simple end, which we denote $\theta(\varphi)$.

The sets $\overline{Q(\varphi)}\cap\mathbb{S}^{k-2}_{\infty}(\bar t)$ form a
singular foliation of $\mathbb{S}^{k-2}_{\infty}(\bar t)$ by $(k-3)-$spheres for $0<\varphi<\pi$ with two singular points
$\theta=\theta(0)$ and $\widetilde\theta =\theta(\pi)$. Thus, each leaf
$\overline{Q(\varphi)}\cap\mathbb{S}^{k-2}_{\infty}(\bar t)$ contains a single
accumulation point $\theta(\varphi)$, the simple end if $0<\varphi<\pi$,
and the points $\theta$ and $\widetilde\theta$ for $\varphi=0$ or $\pi$, when the set
$\overline{Q(\varphi)}\cap\mathbb{S}^{k-2}_{\infty}(\bar t)$
consists of a single point.

\begin{figure}[htbp]
\centering
\includegraphics[width=5.5cm]{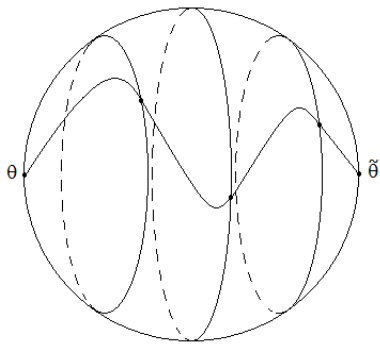}
\caption{The curve of accumulation points in
$\partial P_{\gamma}(\overline{t})$, from $\theta$ to $\widetilde{\theta}$.} \label{thetatil}
\end{figure}

The set formed by these points is the graph of a function $\theta:[0,\pi] \rightarrow \mathbb{S}^{k-2}_{\infty}(\bar t)$, and it is a closed set in the sphere $\mathbb{S}^{k-2}_{\infty}(\bar t) =
\partial_{\infty}\pi(P_{\gamma}(\overline{t}))$.
By an elementary fact of general topology, it follows
that the function $\theta$ is continuous. Therefore the set of all the accumulation points at infinity of $\pi(\Sigma^k)$ form a
continuous curve in $\mathbb{S}^{k-2}_{\infty}(\bar t)$ with extremities $\theta$ and $\widetilde{\theta}$ (See Figure \ref{thetatil}).

Now taking another point of accumulation at infinity between $\theta$ and $\widetilde{\theta}$, say $\theta'$, we can consider the $2-$plane
$\mu'$ generated by $p_{0}, \theta$ and $\theta'$
and let $\xi'$ be the $(k-1)-$plane orthogonal to $\mu'$
passing through $p_{0}$. Repeating the previous argument for the new
foliation by vertical hyperplanes containing the plane $\xi'$ shows that $\theta$ and this other point $\theta'$ should be extremities of another curve which also consists of
all the points of accumulation at infinity of $\pi(\Sigma^k)$. This is
absurd since $\widetilde{\theta}$ and nearby limit points are excluded from this curve. The contradiction shows that $\theta$ is the only accumulation point at infinity.
As in Case 2, every vertical hyperplane whose projection in
$\Hset^k$ does not have $\theta$ as a limit point must meet $\pi(\Sigma^k)$ in a compact set, so $\Sigma^k$ has a simple end at $\theta$.


\bigskip


\noindent{\bf Case 4.} The intersection $\Sigma^{k-1}(t)$ becomes non-compact for some $\bar t>0$.

\vspace{0.5cm}
First, note that this is the only remaining case to complete the
proof of the Main Theorem.

By transversality the set of positive
$t$'s for which the intersection
of $\Sigma^k$ with $P_{\gamma}(t)$ is compact is open in the
line. It follows that there is a first value of $t$, say
$\overline{t}>0$, such that $ \Sigma^{k-1}(\overline {t})$ is not
compact. Then $\Sigma^{k-1}(\overline{t})$ is the limit of
$\Sigma^{k-1}(t) $ as $t$ approaches $\overline{t}$ from below.

Since $\Sigma^{k-1}(\overline{t})$ is not compact, the inductive
hypothesis implies that $ \Sigma^{k-1}(\overline{t})$ is a $(k-1)-$surface homeomorphic to $\mathbb{R}^{k-1}$ and that it is
a vertical graph over a convex domain in $\mathbb{H}^{k-1}$ or has a
simple end.

\begin{lemma} Let $P$ be a vertical hyperplane that meets $\Sigma^k$,
such that $P$ is not the vertical plane $P_0$ tangent to $\Sigma^k$
at the point $p_0$. Then the intersection $\Sigma^k\cap P$ is not a
vertical graph over a convex domain. \label{notvertical}
\end{lemma}

\noindent{\bf Proof.} First, suppose that $p_0\in P$. Then the
tangent plane to $\Sigma^k\cap P$ at the point $p_0$ is a vertical
plane, so $\Sigma^k\cap P$
is not a vertical graph.

Now suppose that $p_0$ is not in $P$ and consider
the complete geodesic $\bar\gamma$ that passes through
$p_0$ and is orthogonal to $P$. For every $t\in\Rset$, let
$P_{\bar\gamma}(t)$ be the vertical hyperplane orthogonal to
$\bar\gamma$ passing through the point $\bar\gamma(t)$,
where $\bar\gamma$ is parametrized so that
$\bar\gamma(0)=p_0$ and $\bar\gamma(t)$ for $t>0$ is on the external side of $P_0$. The hyperplanes
$P_{\bar\gamma}(t)$ form a foliation of
$\HkxR$. For some $t_0<0$, $P_{\bar\gamma}(t_0)$ will be tangent to
$\Sigma^k$, and for $t<t_0$, the intersection $\Sigma^k \cap
P_{\bar\gamma}(t)$ will be empty.

By repeating the arguments of Cases 1, 2, and 3, we see that
the only possibility for the connected set
$\Sigma^k\cap P_{\bar\gamma}(t)$
to be a vertical graph occurs in the situation of Case 4,
that is, when there is a number $t'>0$ such that
$\Sigma^k$ is transverse to $P_{\bar\gamma}(t)$ for $0<t\leq t'$ with a compact intersection for $0<t<t'$,
but $\Sigma^k\cap P_{\bar\gamma}(t')$ is not compact.

In this situation, let $\bar\Sigma^{k-1}(t')$ denote $\Sigma^k\cap P_{\bar\gamma}(t')$ and suppose that
is a vertical graph of a function $f$ over a convex domain
in $\pi(P_{\bar\gamma}(t'))$, which is isometric to $\mathbb{H}^{k-1}$
with $k-1\geq 2$,
in order to obtain a contradiction. Let $\beta$ be a geodesic segment contained in the domain of $f$. For each point $q$ of
$\beta$, the vertical line $r(q)$ passing through the point $(q,0)$
meets $\bar\Sigma^{k-1}(t')$ in a unique point
$(q,f(q))$, and these points form a curve $\widetilde\beta$ lying
over $\beta$, so that $\pi(\widetilde\beta)=\beta$. If $\mu_q$ is
the complete geodesic in $\mathbb{H}^k$
containing the points $\pi(p_0)$ and $q\in\beta$, then the intersection of $\Sigma^k$ with the $2-$plane $\mu_q\times{\mathbb R}$ is a complete curve immersed in
$\mu_q\times{\mathbb R}$ with strictly positive
curvature and a vertical tangent at $p_0$, and
$\mu_q\times{\mathbb R}$ is isometric to ${\mathbb R}^2$.
Since $r(q)$ meets $\bar\Sigma^{k-1}(t')$ in a single point
$(q,f(q))$, at that point two branches of the curve must cross
each other. This holds for every point $q$ in $\beta$, so the curve
$\widetilde\beta$ over $\beta$ must have strictly positive curvature
both above and below in the Euclidean $2-$plane strip
$\beta\times{\mathbb R}$, but that is absurd.
\qed

\bigskip

The lemma shows that $\Sigma^{k-1}(\overline{t})$ cannot be
a vertical graph, so by the inductive hypothesis, it
must have a simple end, which we denote $\theta$,
in $\mathbb{S}^{k-2}_{\infty}=\partial \pi(P_{\gamma}(\overline t))$. We shall show that $\theta$ is also a simple end of $\Sigma^k$.

Let $\Omega(\gamma,\theta)$ be the hyperbolic $2-$plane in
$\mathbb{H}^{k} \times \{0\} $ that contains the complete geodesic
$\gamma$ orthogonal to $\Sigma^k$ at $p_0$ and also has
$(\theta,0)$ as an accumulation point at infinity.
Parametrize the circle $\mathbb{S}^{1}_{\infty} =
\partial_{\infty}(\Omega(\gamma, \theta)) $ by the numbers $0$ to $2\pi$ such that $0$ is the parameter for $(\theta,0)$
and $\pi$ is the other point in $\mathbb{S}^{1}_{\infty}\cap \partial P_{\gamma}(\overline {t})$, with the orientation such that the points $ 0 < s < \pi $ are on the same side
of $\partial P_{\gamma}(\overline {t})$ as the point $ p_{0} $.

Now for $\delta > 0$ near to $0$
consider the complete geodesic $\{\delta, s \}$ from $\delta$ to $s$
in $\Omega(\gamma,\theta)$, where
$\delta < s < 2 \pi $, and let $ W(\delta, s) $ be the vertical $k-$plane in $\mathbb{H}^k \times \mathbb{R}$ that contains the geodesic $\{\delta, s \}$ and is orthogonal to the plane $\Omega(\gamma, \theta)$.
Let us analyze the
intersection $\Sigma^{k-1}(\delta,s)$ of $\Sigma^k$ and $W(\delta,s)$ for $\pi < s < 2\pi$. There are two possible situations:

\begin{figure}[htbp]
\centering
\includegraphics[width=5.0cm]{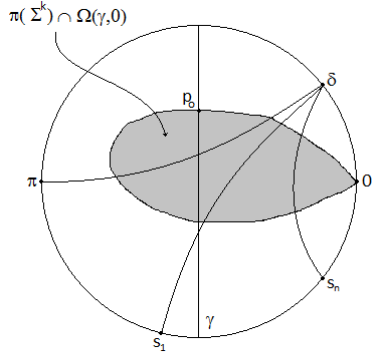}
\caption{Situation 1 of Case 4.}\label{situation1}
\end{figure}

\begin{enumerate}
   \item The hypersurfaces $\Sigma^{k-1}(\delta, s)$ are always compact for all $\pi \leq s < 2\pi$ and all $\delta > 0$ sufficiently near to zero.
   In this case, as $\delta \rightarrow 0$ we have a situation similar to Case 2, so the argument there shows that $\theta$ is a simple end of $\Sigma^k$ (See Figure \ref{situation1}).
 \item There are numbers $\delta >0$ arbitrarily near to zero such that $\Sigma^{k-1}(\delta, s)$ becomes non-compact for some $s, \pi < s < 2\pi $.
\end{enumerate}

\begin{figure}[htbp]
\centering
\includegraphics[width=5.0cm]{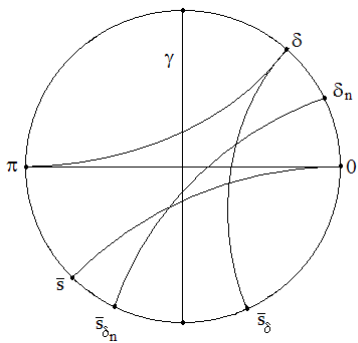}
\caption{Situation 2 of Case 4.}\label{situation2}
\end{figure}

To handle this second situation, for a fixed $\delta>0$ let $\bar{s}_{\delta}$ be the smallest $s$ such that $\Sigma^{k-1}(\delta, \overline{s}_{\delta})$ is not compact. It is easy to see that as $\delta$ decreases $\overline{s}_{\delta}$ cannot increase. Note that there is a number $s$ between $\delta$ and $\pi$ where $ W(\delta, s) $ is tangent to $\Sigma^k$ and at this point $W(\delta, s)$ is a tangent vertical hyperplane. Thus we can apply Lemma \ref{notvertical} and conclude that $ \Sigma^{k-1}(\delta, \overline{s}_{\delta})$ cannot be a vertical graph so it must be a $(k-1)-$surface with a simple end, which we denote
by $\widetilde{\phi}_{\delta}$. Note that $\widetilde{\phi}_{\delta}$ is not necessarily in the circle $\partial_{\infty}(\Omega(\gamma, 0))$, but it is in $S^{k-1}_{\infty}$.
Let $\overline{s}$ be the limit of  $ \overline{s}_{\delta}$
as $\delta$ decreases to $0$. Then there is a decreasing sequence
$\delta_n$ such that $\widetilde{\phi}_{\delta_{n}}$
converges to a point
$ \overline{\phi} \in  S^{k-1}_{\infty} $.
Note that $\overline{s} \geq \pi $ and suppose $ \overline{s} < 2\pi $ to arrive at a contradiction. Let $T_{n}$ be the hyperbolic $(k-1)-$plane that contains the geodesic $ \{ \delta_{n}, \overline{s}_{\delta_{n}} \} $ and is orthogonal to the $2-$plane $\Omega(\gamma, \theta)$. Then the limit $T$ of the $(k-1)-$planes $T_{n}$ is the
hyperbolic $(k-1)-$plane that contains the geodesic
$\{0,\overline{s}\}$ and is orthogonal to $\Omega(\gamma,\theta)$. Note that $ \widetilde{\phi}_{\delta_{n}} \in \partial_{\infty}T_{n} $ and $ \overline{\phi} \in \partial_{\infty}T $. To see that $\overline{\phi} $ is an accumulation point of $ \Sigma^k \cap (T \times \mathbb{R}) $, we use the following proposition.

\begin{prop}
Supposing the hypotheses of the Main Theorem and Hypothesis
\ref{hyp2}, let $\theta\in \partial\Hset^k$
be an accumulation point both of $\pi(\Sigma^k)$ and of a totally
geodesic hyperplane $Q$ in $\Hset^k$.
If there exists a neighborhood $V$ of
$\theta$ in $\bar{\Hset}^k =
\Hset^k\cup\partial {\Hset}^k$ such that
$V\cap\pi(\Sigma^k)\cap Q$ is empty, then $\theta$ is the only
accumulation point of $\pi(\Sigma^k)$ in the sphere at infinity
$\partial {\Hset}^k$. \label{separationprop}
\end{prop}

Supposing this Proposition for the moment, note that both
$\theta$ and $\overline{\phi}$ are accumulation points at
infinity of $T$ and also of $\pi(\Sigma^k)$. If
$\theta$ and $\overline{\phi}$ are distinct, then the
Proposition shows that both of them
must be accumulation points of $\pi(\Sigma^k)\cap T$.
This is impossible since by the inductive hypothesis
$\pi(\Sigma^k)\cap T$ must have a simple end and thus it has
only one accumulation point at infinity. Consequently
$\theta$ and $\overline{\phi}$ must coincide, so $\pi(\Sigma^k)$
has only one accumulation point $\theta$ at infinity.
As in Case 2, $\theta$ must be a simple end of
$\pi(\Sigma^k)$. This will complete the proof of the Main Theorem,
once we have proven the Proposition. \qed

\bigskip

We shall use the following lemma in the proof of Proposition
\ref{separationprop}.

\begin{lemma}
In the conditions of Proposition \ref{separationprop}, the image
$\pi(\Sigma^k)$
is a closed set in $\Hset^k$ and its frontier $\partial\pi(\Sigma^k)$
is a complete smooth hypersurface embedded in $\Hset^k$ with strictly
positive curvature at every point.
\end{lemma}

\noindent{\bf Proof.} First, we show that $\pi(\Sigma^k)$ is closed in
$\Hset^k$. Let $q_0 = \pi(p_0)\in \Hset^k$ be the projection of the
point $p_0$ where we are assuming there is a vertical tangent plane and
consider a geodesic $\beta$ in $\Hset^k$ beginning at $q_0= \beta(0)$,
with $s\geq 0$. If there exists some $s>0$ such that
$\beta(s)\notin\pi(\Sigma^k)$, let
$\bar s>0$ be the largest number such that $z=\beta(s)\in \pi(\Sigma^k)$
for all $s\in [0,\bar s)$. Consider the intersection of the plane
$\beta\times{\mathbb R}$ with $\Sigma^k$. This is a complete convex
embedded curve with a vertical point at $p_0$, and since it does
not diverge to infinity, it must contain exactly one other point with
a vertical
tangent plane, say $p_\beta$, such that the geodesic $\beta$ from
$q_0=\beta(0)$ to $q_\beta=\beta(\bar s)=\pi(p_\beta)$ is the
intersection of $\beta$ with
$\pi(\Sigma^k)$. Let $\widetilde B$ be the set of all these points
$p_\beta$ for all such geodesics $\beta$ plus the point $p_0$.
Then $\widetilde B$ is the set of points in $\Sigma^k$ with vertical
tangent planes and we see that the frontier of $\pi(\Sigma^k)$
is contained in $B=\pi(\widetilde B)\subset\pi(\Sigma^k)$,
so $\pi(\Sigma^k)$ is a closed set in $\Hset^k$ with boundary $B$.

Next we claim that $\widetilde B$ is a smooth hypersurface in $\Sigma^k$.
In a neighborhood $V\subset \Sigma^k$ of the point $p_\beta$ let $N$
be the normal vector field to $V$. Let $\pi':\Hset^k\times
{\mathbb R} \to {\mathbb R}$ be the projection on the second factor, and
let $d\pi'_p: T\Sigma^k_p\to T{\mathbb R}_{f(p)}\equiv {\mathbb R}$ be
its differential at $p\in V$. In view of the strictly
positive curvature of
$\Sigma^k$, $d\pi'_p$ is surjective for every $p\in \widetilde B\cap V$,
so the implicit function theorem shows that near $p$, $\widetilde B\cap V$
is a smooth hypersurface in $V$. Furthermore, $\pi|\widetilde B: \widetilde
B\to B$ is a
diffeomorphism and $B$ is a smooth
hypersurface in $\Hset^k$, as claimed.

It remains to show that $B$ has strictly positive curvature in
$\Hset^k$.
Let $N$ be the external unit normal field to $B$. Since
$\Sigma^k$ has a vertical tangent plane at every point of
$\widetilde B$, the external unit normal
field $\widetilde N$
to $\Sigma^k$ is horizontal along $\widetilde B$ and
coincides with $N$, i.e., at $p\in \widetilde B$,
$d\pi_p\widetilde N(p) =
N(\pi(p))$.
For a tangent vector $\widetilde X\in T_p \widetilde B$ with
$d\pi_p\widetilde X = X$
we have
$$A(X) = -\nabla_{X} =
-\widetilde\nabla_{\widetilde X} = \widetilde A(\widetilde X)$$
where $\nabla$ and $\widetilde\nabla$ are the Riemannian connections
on $\Hset^k$ and $\Hset^k\times{\mathbb R}$ and
$A$ and $\widetilde A$ are the shape operators associated to the
second fundamental forms of $B$ and $\Sigma^k$. Thus the principal
curvatures
of $B$ at the point $\pi(p)$, which are the eigenvalues of $A$, coincide with
the eigenvalues of the restriction of $\widetilde A$ to $\widetilde B$,
the principal curvatures of $\Sigma^k$ along $\widetilde B$ at $p$, and these
are all positive since $\Sigma^k$ has strictly positive curvature.
Thus $B$ has strictly positive curvature.
\qed

\bigskip

\noindent{\bf Proof of Proposition \ref{separationprop}.} By the
lemma, the frontier $B$ of $\pi(\Sigma^k)$ is a smooth
hypersurface in $\Hset^k$ contained in $\pi(\Sigma^k)$
and it has strictly positive curvature.
By hypothesis, near
$\theta$ the image $\pi(\Sigma^k)$ and its frontier $B$ must
are disjoint from the totally geodesic hyperplane $Q$. For any totally
geodesic $2-$plane $\eta$ in $\Hset^k$ that is orthogonal to $Q$
and has $\theta$ as a limit point,
the intersection $\eta\cap Q$ is a geodesic in $\Hset^k$, while $\eta\cap B$
is a strictly convex curve (unless it is empty),
so as the points
on this curve move away from $\theta$ they must become more distant
from $\eta\cap Q$. Outside a tiny neighborhood of the sphere at
infinity, the curvature of the boundary $B$ of $\pi(\Sigma^k)$ is
greater than some positive
number $\epsilon$. This forces $\pi(\Sigma^k)$ to move away from $Q$
uniformly as points in it move away from $\theta$. Hence there
cannot be any other accumulation point of $\pi(\Sigma^k)$ in
addition to $\theta$.
Since $\pi(\Sigma^k)$ is connected, it must lie
entirely on one side of $Q$.
\qed


\section{Examples of simple ends}\label{example}

We use the upper half-space model of hyperbolic $n-$space $\Hset^n$,
$$\Rset^n_+=\{x=(x_1,\dots,x_n)\in \Rset^n\ |\ x_n>0\}$$
with the Riemannian
metric given explicitly as $ds^2= x_n^{-2}\sum_{i=1}^n dx_i^2$.
Then $\Hset^n=\Rset^n_+$ with this metric has constant curvature
$-1$.
We shall give examples of complete embedded hypersurfaces with
strictly positive curvature and with a simple end by explicit computations. The hypersurfaces which we consider
are invariant under the group of parabolic isometries that
leave invariant each horizontal horosphere in $\Rset^n_+$
and fix $\infty$, the point at infinity.
Explicitly, for $a=(a_1,\dots,a_{n-1})\in \Rset^{n-1}$,
consider the parabolic isometry
$f_a:\Rset^n_+\to \Rset^n_+$, defined by setting $f_a(x)=x+(a,0)$,
where $(a,0)=(a_1,\dots,a_{n-1},0)$. Then $f_a$
preserves the horizontal horospheres $\{x_n=b\}$, and
$$F_a:\mathbb{H}^n\times \mathbb{R}\to
\mathbb{H}^n\times \mathbb{R},\ F_a(x,t)=(f_a(x),t)$$
is an isometry for every $a\in\Rset^{n-1}$.

Consider an interval  $\mathcal{I} = (t_{1}, t_{2})$, $-\infty < t_{1}< t_{2} < + \infty$, and a positive smooth function $u:\mathcal{I} \rightarrow
\mathbb{R}$,
\begin{equation}u(t)= c_{1}\ln(t-t_{1}) + c_{2}\ln(t_{2} - t)\label{defu}\end{equation}
where $c_1$ and $c_2$ are negative constants and
\begin{equation}t_2-t_1\leq e^{-1}.
\label{t1t2}\end{equation}
It is
clear that $u(t)$ is positive and $\lim_{t\to t_1}u(t) =
\lim_{t\to t_2}u(t) = +\infty$.
Let $$\alpha=\{(0, 0, ..., 0, u(t), t)\ |\ t\in\mathbb{R}\}
\subset\mathbb{H}^n\times \mathbb{R}$$ be the graph of $u$
as a curve in the $x_{n}t-$plane and consider the set
$$\Sigma^n=\cup_{a\in \Rset^{n-1}} F_a(\alpha)\subset \mathbb{H}^n\times \mathbb{R},$$
parametrized by
\begin{equation}
\varphi(x_{1}, x_{2}, ..., x_{n-1}, t) = (x_{1},
x_{2}, ..., x_{n-1}, u(t), t).
\label{parametrization}
\end{equation}
Clearly $\Sigma^n$
is a complete smooth embedded hypersurface.
We shall show that (\ref{t1t2}) implies that
it has strictly positive curvature and has a simple end.

The hypersurface $\Sigma^n$ is preserved by
$A\times {\rm id}_\mathbb{R}: \mathbb{H}^n\times \mathbb{R}\to\mathbb{H}^n\times \mathbb{R}$ for
every Euclidean isometry $A:\Rset^n_+\to \Rset^n_+$
that preserves $x_n$; this includes not only horizontal translations, but also rotations of $\Rset^n_+$ and compositions, and these all produce isometries of $\mathbb{H}^n\times \mathbb{R}$.

Let $\{\partial_1=\partial/\partial x_1,\dots,
\partial_n=\partial/\partial x_n,
\partial_t=\partial/\partial t\}$ be the usual frame
in $\mathbb{H}^n \times \mathbb{R}$ and let
$\lambda = x_n^{-1}$ be the conformal factor of the metric.
Then
the Riemannian connection $\overline\nabla$ in $\mathbb{H}^n \times \mathbb{R}$ is given by
\begin{displaymath}\begin{array}{lll}
\overline\nabla_{\partial_i}\partial_j =&\delta_{ij}\lambda\partial_n\ &{\rm for} \ i,j< n \nonumber \\
\overline\nabla_{\partial_i}\partial_n =& \overline\nabla_{\partial_n}\partial_i =-\lambda\partial_i &{\rm for}\ i<n \nonumber \\
\overline\nabla_{\partial_n}\partial_n =& -\lambda\partial_n \nonumber
\\ \overline\nabla_{\partial_t}\partial_i\ =&
\overline\nabla_{\partial_i}\partial_t=
\overline\nabla_{\partial_t}\partial_t=0 &{\rm for}\ i\leq n,
\end{array}
\end{displaymath} where  $\delta_{ij}$ is the Kronecker delta.

Using the parametrization (\ref{parametrization}) and
setting $u_t=\partial u/\partial t$ as usual,
we see that the vector fields
\begin{align}
&\varphi_i = \varphi_{x_{i}}  = \partial_i,\ \hspace{2.3cm}{\rm for}\  i=1,\dots,n-1,\ {\rm and}\nonumber \\
&\varphi_n = \varphi_{t} = u_{t}\partial_{n} + \partial_t \nonumber
\end{align}
form a basis of the tangent space to $\Sigma^n$.
We can calculate the coefficients of the first fundamental form $g_{ij} = \langle
\varphi_{i}, \varphi_{j}\rangle $ of $\Sigma^n$ to be
\begin{displaymath}
g_{ij} = \left\{ \begin{array}{ll}
\lambda^2 &{\rm if}\  i=j\ {\rm with} \ i = 1, ..., n-1 \nonumber \\
\lambda^2u_{t}^{2} + 1 &{\rm for} \hspace{0.3cm}i=j=n \nonumber \\
 0 &{\rm for}\ i \neq j\  {\rm with}  \  i, j = 1, ..., n. \nonumber
\end{array}\right.
\end{displaymath}
The unit normal vector field to $\Sigma^n$ is
\begin{eqnarray}
N = (\lambda^{-1}\partial_{n} -
\lambda u_{t}\partial_t)/m  \nonumber
\end{eqnarray}
where $m = \sqrt{1 + \lambda^2 u_{t}^2}$.
Then we compute the coefficients  $b_{ij} = \langle \overline{\nabla}_{\varphi_{j}}\varphi_{i} ,
 N \rangle $ of the second fundamental form to be
\begin{displaymath}
b_{ij} = \left\{ \begin{array}{ll}
\lambda^2/m \hspace{0.3cm} &{\rm if} \hspace{0.3cm}i=j \hspace{0.3cm} {\rm with} \hspace{0.3cm}i = 1, ..., n-1 \nonumber \\
\lambda(-\lambda u_{t}^2 + u_{tt})/m\hspace{0.3cm} &{\rm for} \hspace{0.3cm} i=j=n \nonumber \\
 0  \hspace{0.3cm} &{\rm for} \hspace{0.3cm} i \neq j\hspace{0.3cm} {\rm with} \hspace{0.3cm}i, j = 1, ..., n. \nonumber
\end{array}\right.
\end{displaymath}
Therefore the matrix $(b_{ij})_{n \times n}$ of
coefficients of the second fundamental form is a diagonal matrix with the
eigenvalues
\begin{align}
 &\mu_{i} = \lambda^2/m\ \hspace{2.7cm}{\rm for}\ i=1,\dots, n-1\
 {\rm and}\nonumber \\
 &\mu_{n} = \lambda(-\lambda u_{t}^2  + u_{tt})/m\nonumber
\end{align}
along the diagonal.
Since $m = \sqrt{1 + \lambda^2 u_{t}^2} > 0$, it is clear that $\mu_i>0$ for
$i=1,\dots,n-1$.

\begin{lemma} Condition (\ref{t1t2}) implies that $-u_{t}^2 + uu_{tt}> 0$.
\label{utt}
\end{lemma}
Since $\lambda=1/u>0$ this lemma implies that
$\mu_n$ is also positive. Thus $\Sigma^n$, as defined by the explicit
formulas given above,
is a complete embedded hypersurface in $\mathbb{H}^n \times \mathbb{R}$
with a simple end and
strictly positive curvature.

\bigskip

\noindent{\bf Proof of Lemma \ref{utt}.} Calculate the derivatives
$u_t$ and $u_{tt}$ from  (\ref{defu}), the definition
of $u(t)$. Then a short calculation gives
$$-\lambda u_{t}^2  + u_{tt}=-\frac{c_2^2}{(t-t_1)^2}
[1+\ln(t-t_1)]
-\frac{c_1^2}{(t_2-t)^2} [1+\ln(t_2-t)]$$
$$\qquad\qquad\qquad+\frac{c_1c_2}{(t-t_1)(t_2-t)}\left[2-\frac{(t-t_1)}{(t_2-t)}\ln(t-t_1)
-\frac{(t_2-t)}{(t-t_1)}\ln(t_2-t)\right].$$
Since $t_2-t_1\leq e^{-1}$, $t_1<t<t_2$, and the parameters $c_1$ and $c_2$ are negative,
it is easy to check that each of the terms on the right hand side is
positive. \qed

\end{document}